\documentclass[12pt,reqno]{amsart}
\usepackage[T1]{fontenc}
\usepackage{lmodern}
\usepackage[letterpaper,textwidth=5.90in,textheight=8.40in,centering]{geometry}
\usepackage{amsmath,amssymb,mathtools,mathrsfs}
\usepackage{microtype}
\usepackage{needspace}
\usepackage{tikz}
\usetikzlibrary{arrows.meta}
\usepackage{xcolor}
\definecolor{linkblue}{RGB}{32,63,103}
\usepackage[colorlinks=true,
            linkcolor=blue,
            citecolor=blue,
            urlcolor=blue,
            bookmarksnumbered=true]{hyperref}
\hypersetup{
   unicode=false,          
   pdftoolbar=true,        
bookmarksnumbered=true,
   pdfmenubar=true,        
   pdffitwindow=false,     
   pdfstartview={},    
   pdftitle={},    
   pdfauthor={},     
   pdfsubject={},   
   pdfcreator={},   
   pdfproducer={}, 
   pdfkeywords={}, 
   pdfnewwindow=true,      
   colorlinks=true,       
   linkcolor=blue,          
   citecolor=blue,        
   filecolor=blue,      
   urlcolor=blue          
}

\numberwithin{equation}{section}
\newtheorem{theorem}{Theorem}[section]
\newtheorem{proposition}[theorem]{Proposition}
\newtheorem{lemma}[theorem]{Lemma}

\theoremstyle{definition}

\theoremstyle{remark}
\newtheorem*{remark}{Remark}
\newcommand{\T}{\mathcal T}
\newcommand{\M}{\mathcal M}
\newcommand{\MCG}{\operatorname{Mod}}
\newcommand{\R}{\mathbb R}
\newcommand{\Z}{\mathbb Z}
\newcommand{\HH}{\mathbb H}
\newcommand{\WP}{\mathrm{WP}}
\newcommand{\dwp}{d_{\WP}}
\newcommand{\gwp}{g_{\WP}}
\newcommand{\Ric}{\operatorname{Ric}}
\newcommand{\Hess}{\operatorname{Hess}}
\newcommand{\Vol}{\operatorname{Vol}_{\WP}}
\newcommand{\sech}{\operatorname{sech}}
\newcommand{\inj}{\operatorname{inj}}
\newcommand{\diam}{\operatorname{diam}_{\WP}}
\newcommand{\grad}{\nabla}
\newcommand{\ip}[2]{\left\langle #1,#2\right\rangle}
\newcommand{\norm}[1]{\left\|#1\right\|}
\newcommand{\dd}{\,dt}
\newcommand{\dA}{\,dA}
\newcommand{\eps}{\varepsilon}

\newcommand{\Int}{\operatorname{Int}}
\newcommand{\Rea}{\operatorname{Re}}
\newcommand{\Ima}{\operatorname{Im}}
\title[Weil--Petersson diameter]{The Weil--Petersson diameter of moduli space for large genus}

\author{Yunhui Wu}

\address{The Yau Mathematical Sciences Center\\
\& Department of Mathematical Sciences, \\
Tsinghua University, Beijing 100084, P. R. China}

\email{yunhui\_wu@mail.tsinghua.edu.cn}

\date{\today}
\subjclass[2020]{32G15, 30F60, 53C22}
\keywords{Weil--Petersson metric, moduli space, diameter, large genus}

\begin{document}
\begin{abstract}
We determine the large-genus order of the Weil--Petersson diameter of
$\M_{g,n}$ for each fixed number of punctures: namely,
$\diam(\M_{g,n})\asymp\sqrt g\log g$. The main new ingredient is an
integrated Ricci estimate for geodesics that minimize in moduli space and
have thick endpoints, with no systole bound along their interiors. 
\end{abstract}
\maketitle

\section{Introduction}

Let $S_{g,n}$ be an oriented surface of genus $g$ with $n$ punctures, and
let $\T_{g,n}$ be its Teichm\"uller space equipped with the
Weil--Petersson metric. The mapping class group $\MCG(S_{g,n})$ acts by
isometries, and the quotient
\[
 \M_{g,n}=\T_{g,n}/\MCG(S_{g,n})
\]
is the moduli space of hyperbolic surfaces of type $(g,n)$. We write
$\T_g=\T_{g,0}$ and $\M_g=\M_{g,0}$. The completion of $\T_{g,n}$ is the augmented Teichm\"uller space; its boundary strata parametrize nodal surfaces. The quotient of this completion by the mapping class group is compact. In particular, the diameter $\diam(\mathcal M_{g,n})$ of
$\M_{g,n}$ is finite, despite the incompleteness of its metric. One may see \emph{e.g.}~\cite{W08} for more details.

Our main result is as follows.

\begin{theorem}\label{thm:main}
There exists a universal constant $c>0$ such that, for every fixed
integer $n\geq0$,
\[
   \liminf_{g\to\infty}~
   \frac{\diam(\mathcal M_{g,n})}{\sqrt g\log g}\geq c.
\]
\end{theorem}

\begin{remark}
With normalization \eqref{eq:wp-normalization}, one may take the constant $c$ above to be
\[
   c=\frac1{454\sqrt6}.
\] 
This Theorem \ref{thm:main} together with the work \cite{CP} of Cavendish--Parlier implies that for every fixed $n\geq0$,
\begin{equation}
\diam(\mathcal \M_{g,n})
   \asymp
   \sqrt g\log g. \nonumber
\end{equation}
\end{remark}

The study of asymptotic geometry of $\M_{g,n}$ as $g$ tends to infinity, has an active history. For example, the asymptotic behavior of the Weil--Petersson volume of $\M_{g,n}$ was studied \emph{e.g.} by Penner~\cite{Penner}, Schumacher-Trapani~\cite{ST01}, and Mirzakhani--Zograf \cite{MZ15}. As introduced above, Cavendish-Parlier \cite{CP} studied the asymptotic behavior of $\diam(\M_{g,n})$ and showed that the ratio $\frac{\diam(\M_{g,n})}{\sqrt{g}}$ is bounded below by a positive constant and above by a constant multiple of $\log{g}$. For the upper bound, they refined Brock's quasi-isometry of $\T_{g,n}$ to the pants graph \cite{Brock03}. Brock--Bromberg \cite{BrockBromberg16} showed that the length of a shortest Weil--Petersson closed geodesic in $\M_{g}$ is comparable to $\frac{1}{\sqrt{g}}$. Together with M. Wolf \cite{WW18}, we studied the $\ell^p$-norm $(1\leq p \leq \infty)$ of the Weil--Petersson curvature operator and sectional curvature of $\M_g$ for large genus. In \cite{Wu19}  it was shown that the inradius of $\M_{g,n}$ is comparable to $\sqrt{\log g}$; see also \cite{BridgemanBromberg24,Wu24} for related results. Joint with M. Bridgeman \cite{BW}, we showed that the average total Weil--Petersson scalar curvature is comparable to $-g$; one may also see \cite{ST24} for refined results.

\subsection*{Strategy on the proof of Theorem \ref{thm:main}}
We first prove Theorem \ref{thm:main} for closed surfaces; a pinching argument and the boundary-stratum comparisons of \cite{CP} then transfer the same leading
coefficient to every fixed number of punctures.

For $\eta>0$, set the $\eta$-thick part as
\[
 \M_g^{\geq\eta}
 =\{X\in\M_g:\ell_{\mathrm{sys}}(X)\geq\eta\},
\]
where $\ell_{\mathrm{sys}}$ denotes the length of a shortest closed geodesic.
Our principal estimate concerns geodesics which minimize after projection
to $\M_g$.

\begin{theorem}\label{thm:integrated-intro}
For every $\eta>0$ there is a constant $A_\eta>0$, independent of $g$,
with the following property. Suppose that
$\gamma:[0,L]\to\T_g$ is a unit-speed Weil--Petersson geodesic,
its projection to $\M_g$ is globally length minimizing, and its
endpoints lie over
$\M_g^{\geq\eta}$. Then
\begin{equation}\label{eq:intro-ricci}
 \int_0^L-\Ric_{\WP}(\dot\gamma(t),\dot\gamma(t))\,dt
 \leq A_\eta L.
\end{equation}
Here $\Ric_{\WP}$ is the real Riemannian Ricci tensor.
One may take $A_\eta=C(1+\eta^{-2})$ with $C>1$ universal.
\end{theorem}

The key feature of Theorem~\ref{thm:integrated-intro} is that thickness
is required only at the endpoints; the interior of the geodesic $\gamma$ may
enter the thin part. We outline the proof strategy later in the
introduction.

To derive the volume estimate below, we use radial Jacobi fields
and the Riccati inequality along minimizing geodesics from a thick
basepoint, up to the cut locus in $\M_g$.
Theorem~\ref{thm:integrated-intro} then controls the polar Jacobian
at each thick endpoint. Integrating over the corresponding minimizing
polar domain yields the desired estimate. Here the ball is defined
using the ambient Weil--Petersson distance on $\M_g$, rather than
the intrinsic distance on the thick part.

\begin{theorem}\label{thm:ball-intro}
Let $A_\eta\geq1$ be as in Theorem~\ref{thm:integrated-intro}. Then there is a universal $C>0$ such
that, for all
$X\in\M_g^{\geq\eta}$ and $R>0$,
\begin{equation}\label{eq:intro-ball}
 \Vol\bigl(B_{\M_g}(X,R)\cap\M_g^{\geq\eta}\bigr)
 \leq \exp\bigl(Cg+\sqrt{6A_\eta g}\,R\bigr).
\end{equation}
\end{theorem}

On the other hand, Mirzakhani's integration formula and volume estimates
\cite{Mir13} give
\begin{equation}\label{eq:intro-volume}
 \Vol(\M_g^{\geq1})\geq\exp(2g\log g-Cg).
\end{equation}
In fact, the $1$-thick part has at least one third of the total volume for
all sufficiently large $g$. If $R>\diam(\M_g)$, then a
ball of radius $R$ centered at a $1$-thick surface contains
$\M_g^{\geq1}$. Comparing \eqref{eq:intro-ball} and
\eqref{eq:intro-volume} gives
\[
 2g\log g-Cg\leq Cg+\sqrt{6A_1g}\,R,
\]
hence applying $R\downarrow \diam(\mathcal M_{g})$ gives
\[
   \liminf_{g\to\infty}~
   \frac{\diam(\mathcal M_{g})}{\sqrt g\log g}\geq \frac{2}{\sqrt{6A_1}}\geq c.
\]
Keeping track of constants yields the value of $c$ stated in Theorem~\ref{thm:main}.

\subsection*{Strategy on the proof of Theorem~\ref{thm:integrated-intro}}
We now outline the proof of Theorem~\ref{thm:integrated-intro}. For a simple
closed curve $\alpha\subset S_g$, put
\[
 r_\alpha=\sqrt{\ell_\alpha},\qquad
 \lambda_\alpha=\grad r_\alpha,
 \qquad
 a_\alpha(v)=\ip{\lambda_\alpha}{v},\quad
 b_\alpha(v)=\ip{J\lambda_\alpha}{v}.
\]
Thus, $a_\alpha$ and $b_\alpha$ are respectively the radial and angular
components associated with the root-length $r_\alpha$. Motivated by joint work with Bridgeman~\cite{BW}, we first bound
the singular part of $-\Ric_{\WP}(v,v)$ by
\[
 C\sum_{\alpha\ \mathrm{short}}
 \frac{a_\alpha(v)^2+b_\alpha(v)^2}{r_\alpha^2}.
\]
Wolpert's second variation formula~\cite{W08} and convexity of
$r_\alpha$ bound the integrated radial contribution by the angular
contribution and a remainder bounded by $CL$. The essential input
here is a summed collar estimate whose constant is independent
of the genus and the number of short curves. It therefore remains
to control the angular contribution.

For this, we deform $\gamma$ by time-dependent twists on
short-curve excursions. Each twist parameter vanishes at entry
and is integral at exit, so the comparison path has the same
endpoints in $\M_g$. Global minimality in $\M_g$ then bounds its
energy below by $L$. We choose the twists to reduce the angular
energy and estimate their metric distortion using the same
summed collar bound. Controlling the integer correction on
excursions at scale $2s$ and comparing energies yields
\[
 \sum_\alpha
 \int_{\{t:r_\alpha(\gamma(t))\leq s\}}
 b_\alpha(\dot\gamma(t))^2\,dt
 \leq Cs^4L
\]
for all sufficiently small $s$, with the allowable range depending
only on $\eta$. Layer-cake integration gives the required bound
for the weighted angular terms $b_\alpha^2/r_\alpha^2$.
Together with the radial estimate above, this proves
\eqref{eq:intro-ricci}.

\subsection*{Plan of the paper}
Section~\ref{sec:prelim} fixes the conventions and recalls the needed
geometric facts. Section~\ref{sec:collars} proves estimates which may be
summed over all short collars without loss in the genus, and
Section~\ref{sec:ricci} turns them into a directional Ricci bound.
Sections~\ref{sec:twists} and \ref{sec:integrals} contain the twist
comparison and the integrated curvature estimate. The two volume
estimates are established in Sections~\ref{sec:balls} and
\ref{sec:diameter}; they are combined in Section~\ref{sec:diameter} to
prove Theorem~\ref{thm:main}.

\subsection*{AI disclosure}  All mathematical ideas were developed by the author. AI tools were used to assist in the preparation and editing of this manuscript. The author takes full responsibility for the final article.

\subsection*{Acknowledgments}
The author has had many fruitful conversations with many
people over the last one and half decades on this Weil-Petersson diameter problem. He especially thanks Martin Bridgeman, Jeffrey Brock, Changjie Chen, Maryam Mirzakhani, Hugo Parlier, Michael Wolf and Scott Wolpert for many inspiring conversations on related topics. The author is partially supported by the National Key R \& D Program of China (2025YFA1017500) and NSFC grants No. 12361141813 and 12425107. 

\section{Notation and geometric background}\label{sec:prelim}

Unless stated otherwise, we work with closed surfaces of
genus $g\geq2$. By a curve we mean a free homotopy class of essential
simple closed curves, as well as its geodesic representative. The letter
$C$ denotes a positive universal constant that may change from line to
line. In particular, it is independent of the genus, the surface, and the
number and lengths of the short curves under consideration. A subscript
records an allowed dependence. Constants used only to establish
convergence or finiteness will be distinguished from quantitative ones.

\subsection{The Weil--Petersson metric}

Let $X$ have hyperbolic metric $\rho(z)|dz|^2$. A tangent vector to
$\T_g$ at $X$ is represented by a harmonic Beltrami differential
\[
 \mu=\frac{\overline{\varphi}}{\rho},
 \qquad \varphi(z)\,dz^2\in Q(X),
\]
where $Q(X)$ is the space of holomorphic quadratic differentials. Our
normalization is
\begin{equation}\label{eq:wp-normalization}
 \ip{\mu}{\nu}_{\WP}
 =\Rea\int_X\mu\overline\nu\dA,
 \qquad \norm{\mu}_{\WP}^2=\int_X|\mu|^2\dA.
\end{equation}
We identify real tangent vectors and their harmonic representatives. The
complex structure $J$ is an isometry, is parallel for the Levi--Civita
connection, and acts by multiplication by $i$ on these representatives.

We use the nonpositive hyperbolic Laplacian $\Delta$. On a closed surface
define
\begin{equation}\label{eq:D}
 D=2(2-\Delta)^{-1}.
\end{equation}
This is a positive self-adjoint operator with $D(1)=1$. The positivity may
be understood either through its Green kernel or through the maximum
principle. The real Ricci tensor has nonpositive values on all tangent
vectors, since the Weil--Petersson sectional curvature is negative.

For a pants decomposition $P$, let $(\ell_\alpha,\tau_\alpha)_{\alpha\in P}$
be Fenchel--Nielsen coordinates, where $\tau_\alpha$ is measured in units
of hyperbolic length. A full Dehn twist changes $\tau_\alpha$ by
$\ell_\alpha$. In the normalization \eqref{eq:wp-normalization},
\begin{equation}\label{eq:FN}
 \omega_{\WP}=\frac12\sum_{\alpha\in P}
 d\ell_\alpha\wedge d\tau_\alpha,
 \qquad
 t_\alpha:=\frac{\partial}{\partial\tau_\alpha}
 =\frac12J\grad\ell_\alpha.
\end{equation}
These formulas are Wolpert's twist--length duality
\cite{W82,W08}. Changing the overall normalization of the metric
only changes universal constants in our conclusions.

\subsection{Minimizers and boundary strata}

Write $\pi:\T_g\to\M_g$ for the quotient map. The Weil--Petersson
completion of $\M_g$ is a compact length space. Its interior is
geodesically convex, and an interior minimizing path lifts to a
Weil--Petersson geodesic after choosing its initial lift; see
\cite[Theorem~4.18 and the subsequent nonrefraction discussion]{W08}
and \cite[Section~5]{CP}. Thus, for $X,Y\in\M_g$ and a prescribed
lift $\widetilde X$, one can choose a lift $\widetilde Y$ such that
\begin{equation}\label{eq:quotient-min}
 d_{\T_g}(\widetilde X,\widetilde Y)=d_{\M_g}(X,Y).
\end{equation}
We use these minimizing lifts throughout. Any competitor with the same
projected endpoints has length at least this common distance.

If $\sigma$ is a multicurve, the distance from $X$ to the closed stratum $\overline{\T_g(\sigma)}$
where all elements of $\sigma$ are pinched satisfies
\begin{equation}\label{eq:pinching}
 \dwp(X,\overline{\T_g(\sigma)})
 \leq\left(2\pi\sum_{\alpha\in\sigma}\ell_\alpha(X)\right)^{1/2}.
\end{equation}
When $\sigma$ is a pants decomposition, the corresponding marked
stratum is a single maximally noded surface.
The estimate follows from Wolpert's length-gradient bound and pinching
construction \cite{W08}; see also \cite[Theorem~2.1]{CP}.

\subsection{Short collars and cyclic covers}

Fix the numerical cutoff $\ell_0=1/10$ throughout the paper.
The standard collar about a geodesic of length $\ell$ has half-width
\begin{equation}\label{eq:collarwidth}
 w(\ell)=\operatorname{arcsinh}\frac1{\sinh(\ell/2)}.
\end{equation}
Standard collars about disjoint geodesics are embedded and pairwise
disjoint \cite[Chapter~4]{Buser}. The collar lemma also implies that
curves of length at most $\ell_0$ are pairwise disjoint.

The cyclic cover associated to a curve $\alpha$ of length $\ell$ is
\[
 A_\alpha=\{w=x+iy:0<y<\pi\}/(w\sim w+\ell),
 \qquad ds^2=\frac{|dw|^2}{\sin^2y}.
\]
The core is $y=\pi/2$. If $d$ is the distance to this core, then
\begin{equation}\label{eq:sech}
 \sech d=\sin y.
\end{equation}
Inside the standard collar $\widehat C_\alpha$ we use the smaller collar
\begin{equation}\label{eq:innercollar}
 C_\alpha=\{2\ell\leq y\leq\pi-2\ell\},
 \qquad W_\alpha=\log\cot\ell.
\end{equation}
Here $W_\alpha>1$ is its half-width. Since
$\sin(2\ell)>\tanh(\ell/2)$ for $0<\ell\leq\ell_0$, we have
$C_\alpha\subset\widehat C_\alpha$. The cyclic injectivity radius is
\begin{equation}\label{eq:inj-cylinder}
 \sinh\inj_{A_\alpha}(x+iy)
 =\frac{\sinh(\ell/2)}{\sin y}.
\end{equation}
On $C_\alpha$, the right side is at most
$\cosh(1/20)/(4\cos(1/5))<1/\sqrt3$. Hence this collar lies in
$\{\inj_{A_\alpha}<\eps_2\}$, where $\eps_2=\frac12\log3$.
All collar choices are therefore independent of the topology.

If $\mu$ is harmonic, its lift to $A_\alpha$ can be written as
\begin{equation}\label{eq:fourier}
 \mu=\overline{f_\alpha(w)}\sin^2y,\qquad
 f_\alpha(w)=\sum_{k\in\Z}c_{\alpha,k}e^{2\pi ikw/\ell}.
\end{equation}
Let $f_\alpha^\perp=f_\alpha-c_{\alpha,0}$. We shall use the quantities
\begin{equation}\label{eq:Qperp}
 Q_\alpha^\perp(\mu)
 =\int_0^\ell\int_0^\pi
 |f_\alpha^\perp(x+iy)|^2\sin^4y\,dy\,dx,
 \qquad
 N_\alpha^2=\int_{\widehat C_\alpha}|\mu|^2\dA.
\end{equation}
The function $f_\alpha^\perp$ belongs to the cyclic cover. It does not,
in general, define a quadratic differential on $X$.

\section{Uniform estimates over short collars}\label{sec:collars}

Throughout this section, $\sigma$ is a collection of disjoint simple closed geodesics with
\begin{equation}
    \ell_\alpha\leq\ell_0=\frac{1}{10}.
\end{equation}
We combine a local bound for the nonprincipal
Fourier modes with an exterior packing estimate to obtain bounds
independent of $|\sigma|$.

\subsection{An exterior packing estimate}

Let $\mathscr L_\alpha$ be the set of lifts of $\alpha$ to $\HH$.
For $L\in\mathscr L_\alpha$, let $C_L$ be the closed strip of half-width
$W_\alpha$ about $L$. Their interiors are pairwise disjoint: they are lifts of the embedded,
pairwise disjoint inner collars.

\begin{lemma}\label{lem:packing}
For every $z\in\HH$,
\begin{equation}\label{eq:packing}
 \sum_{\alpha\in\sigma}\ \sum_{\substack{L\in\mathscr L_\alpha\\z\notin\Int C_L}}
 \frac{\sech^2 d(z,L)}{\ell_\alpha}\leq76.
\end{equation}
\end{lemma}

\begin{proof}
Fix $z\in\HH$. For each axis in the sum, set
$u_L=d(z,L)-W_\alpha\geq0$. Along the perpendicular path from $z$ to $L$,
choose $p_L$ at distance $1/2$ inside the nearer boundary of $C_L$.
Then $d(z,p_L)=u_L+1/2$, and the balls $B(p_L,1/4)$ lie in the
disjoint strip interiors.

Let $\mathscr A_j$ consist of the axes with $j\leq u_L<j+1$.
Their balls lie in $B(z,j+2)$, so comparison of hyperbolic areas gives
\[
 \#\mathscr A_j\leq
 \frac{\cosh(j+2)-1}{\cosh(1/4)-1}.
\]
Since $e^{-W_\alpha}=\tan\ell_\alpha\leq(101/100)\ell_\alpha$,
\[
 \ell_\alpha^{-1}\sech^2d(z,L)
 \leq4\ell_\alpha^{-1}e^{-2W_\alpha-2u_L}
 \leq4(101/100)^2\ell_0e^{-2u_L}.
\]
Summing over the shells and using $\cosh(j+2)-1\leq e^{j+2}/2$
bounds the LHS of Equation \eqref{eq:packing} by
\[
 \frac{4(101/100)^2\ell_0e^2}
 {2(\cosh(1/4)-1)}
 \sum_{j=0}^\infty e^{-j}
 =
 \frac{4(101/100)^2\ell_0e^2}
 {2(\cosh(1/4)-1)(1-e^{-1})}<76
\]
as desired.
\end{proof}

\subsection{The nonprincipal Fourier modes}We use Teo's explicit constant
\cite[Proposition~3.1]{Teo09},
\begin{equation}\label{eq:Teo-constant}
 C(s)
 =\left[\frac{4\pi}{3}
       \bigl(1-\sech^6(s/2)\bigr)\right]^{-1/2},
 \qquad s>0.
\end{equation}
The disk argument in that proof gives the local estimate
\[
 |\mu(z)|^2\leq C(s)^2
 \int_{B(z,s)}|\mu|^2\dA
 \qquad(0<s<\inj_X(z));
\]
see also \cite[(3.1) and Lemma~3.1]{BW}.
Bridgeman--Wu \cite[Proposition~3.3(3)]{BW} bound each of the
positive and negative Fourier parts on the subcollar
$\inj\leq\eps_2$ by $C(\eps_2)N_\alpha$.
Since $\eps_2=\frac12\log3$ and
$C(\eps_2)<2$, their sum satisfies
\begin{equation}\label{eq:BWlocal}
 |f_\alpha^\perp(w)|\sin^2y\leq4N_\alpha
 \qquad(w\in C_\alpha).
\end{equation}
The norm on the right is the norm on the individual standard collar.
This local form is needed when the inequalities are summed.

\begin{proposition}\label{prop:fourier}
For every harmonic Beltrami differential $\mu\in T_X\T_g$,
\begin{equation}\label{eq:fourier-sum}
 \sum_{\alpha\in\sigma}\ell_\alpha|c_{\alpha,0}|^2
 +\sum_{\alpha\in\sigma}\frac{Q_\alpha^\perp(\mu)}{\ell_\alpha}
 \leq128\norm{\mu}_{\WP}^2.
\end{equation}
\end{proposition}

\begin{proof}
We first estimate the principal coefficient. In the cylinder model,
\[
 |\mu|^2\dA=|f_\alpha(x+iy)|^2\sin^2y\,dx\,dy.
\]
Parseval's identity on each horizontal circle gives
\[
 \int_0^{\ell_\alpha}|f_\alpha(x+iy)|^2\,dx
 =\ell_\alpha\sum_{k\in\Z}|c_{\alpha,k}|^2
       e^{-4\pi ky/\ell_\alpha}
 \geq\ell_\alpha|c_{\alpha,0}|^2.
\]
The interval $\pi/4\leq y\leq3\pi/4$ lies in every standard collar
under consideration, and
$\int_{\pi/4}^{3\pi/4}\sin^2y\,dy\geq\pi/4$. Integrating the last
inequality on this interval yields
\begin{equation}\label{eq:principal-collar}
     \ell_\alpha|c_{\alpha,0}|^2\leq\frac4\pi N_\alpha^2.
\end{equation}

\noindent Since the standard collars are disjoint, we have 
\begin{equation}\label{eq:principal-sum}
 \sum_{\alpha\in\sigma}\ell_\alpha|c_{\alpha,0}|^2
 \leq \frac4\pi\sum_{\alpha\in\sigma}N_\alpha^2
 \leq \frac4\pi\norm{\mu}_{\WP}^2.
\end{equation}

Next split the cylinder integral defining $Q_\alpha^\perp$ into the
inner collar $C_\alpha$ and its complement. On $C_\alpha$,
\eqref{eq:BWlocal} says that the integrand
$|f_\alpha^\perp|^2\sin^4y$ is bounded pointwise by
$16N_\alpha^2$. The Euclidean area of the parameter rectangle for
$C_\alpha$ is
$\ell_\alpha(\pi-4\ell_\alpha)<\pi\ell_\alpha$, and hence
\[
 Q_{\alpha,\mathrm{in}}^\perp\leq16\pi\ell_\alpha N_\alpha^2,
\]
which implies that
\begin{equation}\label{eq:nonprincipal-inner}
 \sum_{\alpha\in\sigma}
 \frac{Q_{\alpha,\mathrm{in}}^\perp}{\ell_\alpha}
 \leq16\pi\norm{\mu}_{\WP}^2.
\end{equation}

The exterior part requires a different argument, because a pointwise
estimate on one collar would not remain uniform after summing over
$\sigma$. For every $y$, orthogonality of the Fourier modes gives
\[
 \int_0^{\ell_\alpha}|f_\alpha^\perp(x+iy)|^2\,dx
 \leq\int_0^{\ell_\alpha}|f_\alpha(x+iy)|^2\,dx.
\]
Thus we may replace the nonprincipal part by the full differential
before unfolding. To see the unfolded integrand explicitly, note that
on the cyclic cover
\[
 |\mu|^2\sech^2d(\,\cdot\,,L)\,\dA
 =|f_\alpha|^2\sin^4y\,dx\,dy
\]
for the distinguished lift of the core, by \eqref{eq:sech}. Summing
over the deck translates of a fundamental domain and then over
$\alpha$, with Tonelli's theorem justifying the nonnegative sums, gives
\begin{align}\label{eq:nonprincipal-outer}
 \sum_{\alpha\in\sigma}\frac{Q_{\alpha,\mathrm{out}}^\perp}{\ell_\alpha}
 &\leq\int_X|\mu(z)|^2
 \sum_{\alpha\in\sigma}\ \sum_{\substack{L\in\mathscr L_\alpha\\
 \widetilde z\notin\Int C_L}}
 \frac{\sech^2d(\widetilde z,L)}{\ell_\alpha}\dA(z)
 \leq76\norm{\mu}_{\WP}^2,
\end{align}
where $\widetilde z$ is any lift of $z$; the last inequality is
Lemma~\ref{lem:packing}.

Adding \eqref{eq:principal-sum}, \eqref{eq:nonprincipal-inner}, and
\eqref{eq:nonprincipal-outer} proves the claim.
\end{proof}

The two terms in \eqref{eq:fourier-sum} separate the principal
root-length directions from the transverse Fourier modes. Both estimates
remain uniform after summing over all short collars.

\subsection{The positive Hessian remainder}

For a curve $\alpha$, recall that
\[
 r_\alpha=\sqrt{\ell_\alpha},\qquad
 \lambda_\alpha=\grad r_\alpha,
 \qquad
 a_\alpha(v)=\ip{\lambda_\alpha}{v},\quad
 b_\alpha(v)=\ip{J\lambda_\alpha}{v}.
\]
Define
\begin{equation}\label{eq:Hdef}
 H_\alpha=\Hess r_\alpha-\frac3{r_\alpha}
 (J\lambda_\alpha)^\flat\otimes(J\lambda_\alpha)^\flat,
\end{equation}
where $\flat$ denotes the metric dual. The chain rule gives
\begin{equation}\label{eq:H-Wolpert}
 H_\alpha(v,v)=\frac{
 2\ell_\alpha\Hess\ell_\alpha(v,v)
 -d\ell_\alpha(v)^2-3d\ell_\alpha(Jv)^2}{4r_\alpha^3}.
\end{equation}
By \cite[Corollary~3.4]{W08}, this form is positive semidefinite.
In particular, $\Hess r_\alpha\geq0$, so root length is convex along
WP geodesics. The following estimate controls its remainder uniformly.

\begin{proposition}\label{prop:Hformula}
Let $\ell_\alpha\leq\ell_0$ and let $v\in T_X\T_g$ have the
Fourier expansion \eqref{eq:fourier}. Then
\begin{equation}\label{eq:firstvariation}
 a_\alpha(v)^2+b_\alpha(v)^2
 =\frac{\ell_\alpha}{4}|c_{\alpha,0}|^2,
\end{equation}
and
\begin{equation}\label{eq:Hbound}
 0\leq\frac{H_\alpha(v,v)}{r_\alpha}
 \leq\frac2\pi\frac{Q_\alpha^\perp(v)}{\ell_\alpha}.
\end{equation}
\end{proposition}

We separate the proof into a Fourier estimate and a variation formula.
Throughout the next two lemmas, integration on
$R_\ell=[0,\ell]\times(0,\pi)$ uses Euclidean area $dx\,dy$.
For a periodic holomorphic function
$f(w)=\sum_k c_ke^{2\pi ikw/\ell}$, set
\[
 Af=\sum_k\frac{c_k}{2\pi ik/\ell-1}e^{2\pi ikw/\ell}.
\]

\begin{lemma}\label{lem:multiplier}
If $f^\perp=\sum_{k\ne0}c_ke^{2\pi ikw/\ell}$ and
$\int_{R_\ell}|f^\perp|^2\sin^4y<\infty$, then
\begin{equation}\label{eq:multiplier}
 \int_{R_\ell}|Af^\perp|^2\sin^2y
 \leq\frac43\int_{R_\ell}|f^\perp|^2\sin^4y.
\end{equation}
\end{lemma}

\begin{proof}
For $\xi\ne0$, put $I_j(\xi)=\int_0^\pi e^{-2\xi y}\sin^j y\,dy$.
Expanding $\sin^2y$ and $\sin^4y$ into cosines gives
\[
 I_2(\xi)=\frac{1-e^{-2\pi\xi}}{4\xi(\xi^2+1)},
 \qquad
 I_4(\xi)=\frac{3(1-e^{-2\pi\xi})}
 {4\xi(\xi^2+1)(\xi^2+4)}.
\]
Thus, for the mode $f_k(w)=c_ke^{2\pi ikw/\ell}$ with
$\xi=2\pi k/\ell$,
\begin{equation}\label{eq:single-mode}
 \int_{R_\ell}|Af_k|^2\sin^2y
 =\frac{\xi^2+4}{3(\xi^2+1)}
   \int_{R_\ell}|f_k|^2\sin^4y
 \leq\frac43\int_{R_\ell}|f_k|^2\sin^4y.
\end{equation}
Orthogonality in $x$ proves \eqref{eq:multiplier} for finite sums. It then follows from
Parseval's identity that the Fourier
partial sums converge in $L^2(\sin^4y\,dx\,dy)$. Applying
\eqref{eq:single-mode} to their differences gives convergence of their
images in $L^2(\sin^2y\,dx\,dy)$, with limit $Af^\perp$.
Passing to the limit proves the claim.
\end{proof}

\begin{lemma}\label{lem:Hvariation}
Under the hypotheses of Proposition~\ref{prop:Hformula},
\begin{equation}\label{eq:principal-variation}
 a_\alpha(v)=\frac{r_\alpha}{2}\Rea c_{\alpha,0},
 \qquad
 b_\alpha(v)=-\frac{r_\alpha}{2}\Ima c_{\alpha,0},
\end{equation}
and
\begin{equation}\label{eq:Hexact}
 H_\alpha(v,v)=\frac1{2\pi r_\alpha}
 \left(2\int_{R_{\ell_\alpha}}|Af_\alpha^\perp|^2\sin^2y
       -\Rea\int_{R_{\ell_\alpha}}(Af_\alpha^\perp)^2\sin^2y\right).
\end{equation}
\end{lemma}

\begin{proof}
Suppress the index $\alpha$. Our holomorphic coefficient is $-4$
times Wolpert's coefficient. Consequently,
\cite[Theorem~3.3]{W08} reads
\[
 d\ell(v)=\ell\Rea c_0,\qquad
 \Hess\ell(v,v)=\frac1\pi
 \left(2\int_{R_\ell}|Af|^2\sin^2y
       -\Rea\int_{R_\ell}(Af)^2\sin^2y\right).
\]
The second variation is the Riemannian Hessian because harmonic
Beltrami coordinates are WP normal \cite[Section~2.3]{W08}.
Since $J\mu=i\mu$, the holomorphic coefficient representing $Jv$
is $-if$. The identities $dr=d\ell/(2r)$ and $dr(Jv)=-b(v)$
therefore give \eqref{eq:principal-variation}.

Write $f=c_0+f^\perp$. Since $A(c_0)=-c_0$, orthogonality in $x$
eliminates the mixed terms. The constant-mode contribution to
$\Hess\ell(v,v)$ is
\[
 \frac{\ell}{2}(\Rea c_0)^2
 +\frac{3\ell}{2}(\Ima c_0)^2
\]
and cancels in
\eqref{eq:H-Wolpert}, leaving \eqref{eq:Hexact}.
All integrals are finite by Proposition~\ref{prop:fourier} and
Lemma~\ref{lem:multiplier}. The same weighted $L^2$ convergence
justifies the bilinear integral, since
$\|u^2-v^2\|_{L^1}\leq\|u-v\|_{L^2}\|u+v\|_{L^2}$.
\end{proof}

\begin{proof}[Proof of Proposition~\ref{prop:Hformula}]
Equation~\eqref{eq:firstvariation} follows from
\eqref{eq:principal-variation}. Set
\[
 U=\int_{R_{\ell_\alpha}}|Af_\alpha^\perp|^2\sin^2y.
\]
The absolute value of the second integral in \eqref{eq:Hexact}
is at most $U$. Hence
\[
 0\leq H_\alpha(v,v)
 \leq\frac{3U}{2\pi r_\alpha}
 \leq\frac{2}{\pi r_\alpha}Q_\alpha^\perp(v),
\]
where the last inequality is Lemma~\ref{lem:multiplier}.
Dividing by $r_\alpha$ proves \eqref{eq:Hbound}.
\end{proof}

\begin{proposition}\label{prop:summed-H}
For every $v\in T_X\T_g$,
\begin{align}
 \sum_{\alpha\in\sigma}
 \bigl(a_\alpha(v)^2+b_\alpha(v)^2\bigr)
 &\leq \norm{v}_{\WP}^2,\label{eq:Bessel}\\
 \sum_{\alpha\in\sigma}\frac{H_\alpha(v,v)}{r_\alpha}
 &\leq81\norm{v}_{\WP}^2.\label{eq:Hsum}
\end{align}
For any real $(e_\alpha)_{\alpha\in\sigma}$,
\begin{equation}\label{eq:dual-Bessel}
 \left\|\sum_{\alpha\in\sigma}e_\alpha J\lambda_\alpha\right\|_{\WP}^2
 \leq\sum_{\alpha\in\sigma}e_\alpha^2.
\end{equation}
\end{proposition}

\begin{proof}
Equations~\eqref{eq:principal-sum} and \eqref{eq:firstvariation} give
\[
 \sum_{\alpha\in\sigma}(a_\alpha^2+b_\alpha^2)
 \leq\frac1\pi\norm{v}_{\WP}^2\leq\norm{v}_{\WP}^2.
\]
Similarly, \eqref{eq:nonprincipal-inner},
\eqref{eq:nonprincipal-outer}, and \eqref{eq:Hbound} give
\[
 \sum_{\alpha\in\sigma}\frac{H_\alpha(v,v)}{r_\alpha}
 \leq\left(32+\frac{152}{\pi}\right)\norm{v}_{\WP}^2
 \leq81\norm{v}_{\WP}^2.
\]
For \eqref{eq:dual-Bessel}, set
$w=\sum_{\alpha\in\sigma}e_\alpha J\lambda_\alpha$.
Cauchy--Schwarz and \eqref{eq:Bessel} give, for every $v$,
\[
 |\ip{w}{v}|
 =\left|\sum_{\alpha\in\sigma}e_\alpha b_\alpha(v)\right|
 \leq\left(\sum_{\alpha\in\sigma}e_\alpha^2\right)^{1/2}
      \norm{v}_{\WP}.
\]
Taking the supremum over unit vectors $v$ proves the claim.
\end{proof}

\section{A directional Ricci bound}\label{sec:ricci}

Fix $0<\eps\leq\ell_0$ and let
\[
 \sigma_\eps(X)=\{\alpha:\ell_\alpha(X)\leq\eps\}.
\]
The collar lemma makes $\sigma_\eps(X)$ a multicurve. The purpose of this
section is to isolate the only part of the negative Ricci tensor which
can become unbounded: its components in the radial and twist directions
associated with $\sigma_\eps(X)$.

\subsection{The Bergman density}

Let $(\nu_j)_{j=1}^{3g-3}$ be a Hermitian orthonormal basis of harmonic
Beltrami differentials and put
\[
 \mathcal B_X(z)=\sum_{j=1}^{3g-3}|\nu_j(z)|^2.
\]
Evaluation at $z$ is a complex linear map to a one-dimensional Hermitian
space, so (see \textit{e.g.}~\cite[Proposition~4.1]{BW})
\begin{equation}\label{eq:Bergman-eval}
 \mathcal B_X(z)=\sup_{\norm{\mu}_{\WP}=1}|\mu(z)|^2.
\end{equation}
Using this identity instead of estimating the basis elements separately
avoids a factor $3g-3$.

For $p>1$, define a function on $X$ by
\begin{equation}\label{eq:kernel}
 K_{\alpha,p}(z)=\sum_{L\in\mathscr L_\alpha}
 \sech^p d(\widetilde z,L).
\end{equation}
The sum is independent of the lift $\widetilde z$. For a fixed surface
and curve it converges locally uniformly with its first two derivatives.
Indeed, disjoint balls in the lifted collar bound the number of axes
at distance at most $R$ by $C_{X,\alpha}e^R$; the summands and their
derivatives are $O(e^{-pR})$. The constant $C_{X,\alpha}$ is used only
for convergence, not in the estimates below.

\begin{lemma}\label{lem:Bergman}
For $0<\eps\leq\ell_0$ and all $z\in X$,
\begin{equation}\label{eq:Bergman-bound}
 \mathcal B_X(z)\leq6\eps^{-2}+
 3\sum_{\alpha\in\sigma_\eps(X)}\ell_\alpha^{-1}K_{\alpha,4}(z).
\end{equation}
\end{lemma}

\begin{proof}
Suppose that $z\in C_\alpha$ for some
$\alpha\in\sigma_\eps(X)$. Write the lift of a unit harmonic
differential as the sum of its principal and nonprincipal terms:
\[
 \mu=\overline{c_{\alpha,0}}\sin^2y
     +\overline{f_\alpha^\perp(w)}\sin^2y.
\]
Since $N_\alpha\leq\norm{\mu}_{\WP}=1$,
\eqref{eq:principal-collar} gives
$|c_{\alpha,0}|\leq2/(\sqrt{\pi\ell_\alpha})$, while
\eqref{eq:BWlocal} bounds the second term by $4$. Hence
\[
 |\mu(z)|\leq4+\frac2{\sqrt{\pi\ell_\alpha}}\sin^2y,
 \qquad
 |\mu(z)|^2\leq32+\frac3{\ell_\alpha}\sin^4y.
\]
For the second inequality, we used $(A+B)^2\leq2A^2+2B^2$ and
$8/\pi<3$. The chosen lift of the core contributes
$\sin^4y$ to $K_{\alpha,4}(z)$ by \eqref{eq:sech}. Since
$32\leq6\eps^{-2}$, taking the supremum in
\eqref{eq:Bergman-eval} proves \eqref{eq:Bergman-bound} at every point
of the chosen inner collars.

If $z$ lies outside these collars, we claim that
\begin{equation}\label{eq:inj-outside-inner}
 \inj_X(z)\geq\eps/2.
\end{equation}
Suppose otherwise. Since $\eps/2\leq1/20<\eps_2$, the thin-part
decomposition places $z$ in the standard collar of a primitive closed
geodesic $\alpha$ of length
$\ell_\alpha \leq 2\inj_X(z)<\eps$. In this collar the injectivity radius
agrees with the cyclic injectivity radius. If $z$ were outside
$C_\alpha$, formula \eqref{eq:inj-cylinder}, together with
$\sin(2\ell_\alpha)\leq2\ell_\alpha$ and
$\sinh(\ell_\alpha/2)\geq\ell_\alpha/2$, would give
\[
 \sinh\inj_X(z)\geq\frac{\sinh(\ell_\alpha/2)}
 {\sin(2\ell_\alpha)}\geq\frac14,
\]
whereas $\inj_X(z)<1/20$ implies
$\sinh\inj_X(z)<\sinh(1/20)<1/4$. This contradiction proves
\eqref{eq:inj-outside-inner}. Apply Teo's local estimate with radius $\eps/4$ to obtain
\[
 |\mu(z)|^2
 \leq C(\eps/4)^2
       \int_{B(z,\eps/4)}|\mu|^2\dA.
\]
For completeness, the numerical bound follows from
$R=\tanh(\eps/8)$ and \eqref{eq:Teo-constant}:
\[
 C(\eps/4)^2
 =\frac{3}{4\pi[1-(1-R^2)^3]}
 \leq\frac1{4\pi R^2(1-R^2)^2}
 <6\eps^{-2}.
\]
Indeed, $\tanh x\geq x/(1+x)$ and $\tanh x\leq x$ give
$R\geq10\eps/81$ and $1-R^2>99/100$ for $\eps\leq1/10$.
Taking the supremum over unit $\mu$ completes the proof.
\end{proof}

\subsection{The resolvent and the curvature tensor}
Recall that $ D=2(2-\Delta)^{-1}$ is a positive self-adjoint operator on $X$.
\begin{lemma}\label{lem:resolvent}
For every curve $\alpha$,
\begin{equation}\label{eq:resolvent}
 D(K_{\alpha,4})=\frac12K_{\alpha,2}.
\end{equation}
Moreover, if $\mu$ has expansion \eqref{eq:fourier}, then
\begin{equation}\label{eq:unfold-Q}
 \int_X|\mu|^2K_{\alpha,2}\dA
 =Q_\alpha^\perp(\mu)+\frac{3\pi}{8}\ell_\alpha|c_{\alpha,0}|^2.
\end{equation}
\end{lemma}

\begin{proof}
In Fermi coordinates about a geodesic in $\HH$, the Laplacian on a
function of signed distance $s$ is $h''(s)+\tanh(s)h'(s)$. Direct
differentiation gives
\[
 (2-\Delta)\sech^2s=4\sech^4s.
\]
Apply this identity to the distance from each axis
$L\in\mathscr L_\alpha$. The local $C^2$ convergence established after
\eqref{eq:kernel} permits termwise differentiation and gives
\[
 (2-\Delta)K_{\alpha,2}=4K_{\alpha,4}.
\]
The operator $2-\Delta$ is positive and invertible on the closed
surface $X$. Applying $2(2-\Delta)^{-1}$ to the last identity proves
\eqref{eq:resolvent}.

For the second assertion, choose a fundamental domain in $\HH$ for the
surface group and unfold the nonnegative integral to the cyclic cover with respect to $\alpha$. Equations
\eqref{eq:sech} and \eqref{eq:fourier} give
\[
 \int_X|\mu|^2K_{\alpha,2}\dA
 =\int_0^{\ell_\alpha}\int_0^\pi|f_\alpha|^2\sin^4y\,dy\,dx.
\]
For each fixed $y$, Parseval's identity separates the zeroth mode from
all other modes. The latter contribute precisely
$Q_\alpha^\perp(\mu)$, while the former contributes
\[
 \ell_\alpha|c_{\alpha,0}|^2
 \int_0^\pi\sin^4y\,dy
 =\frac{3\pi}{8}\ell_\alpha|c_{\alpha,0}|^2.
\]
This is \eqref{eq:unfold-Q}.
\end{proof}

\begin{proposition}\label{prop:Ricci}
For every $0<\eps\leq\ell_0$ and every $v\in T_X\T_g$,
\begin{equation}\label{eq:directional-Ricci}
 -\Ric_{\WP}(v,v)\leq 32\eps^{-2}\norm{v}_{\WP}^2+
30\sum_{\alpha\in\sigma_\eps(X)}
 \frac{a_\alpha(v)^2+b_\alpha(v)^2}{r_\alpha^2}.
\end{equation}
\end{proposition}

\begin{proof}
For $\norm{\mu}_{\WP}=1$, choose the Hermitian orthonormal basis with
$\nu_1=\mu$. Sum the sectional-curvature estimate \cite[(2.1)]{BW} over
$\{\nu_j,J\nu_j\}_{j=1}^{3g-3}\setminus\{\mu\}$, an orthonormal
basis of the real orthogonal complement of $\mu$.
The Ricci tensor is the trace of the curvature operator in this real
basis. Moreover, multiplication by $J$ preserves pointwise norm, so the
sum of the pointwise squared norms of these basis elements is at most
$2\mathcal B_X$. The curvature estimate, the positivity of $D$, and
the self-adjointness of $D$ therefore give
\begin{equation}\label{eq:Ricci-Bergman}
 -\Ric_{\WP}(\mu,\mu)
 \leq4\int_X|\mu|^2D(\mathcal B_X)\dA.
\end{equation}
Homogeneity extends the estimate to arbitrary $\mu$.

Apply $D$ to \eqref{eq:Bergman-bound}. Positivity of $D$ preserves the
inequality, while $D(1)=1$ and Lemma~\ref{lem:resolvent} yield
\[
 D(\mathcal B_X)\leq6\eps^{-2}
 +\frac32\sum_{\alpha\in\sigma_\eps(X)}
 \ell_\alpha^{-1}K_{\alpha,2}.
\]
Substitution in \eqref{eq:Ricci-Bergman} gives
\[
\begin{aligned}
 &   -\Ric_{\WP}(\mu,\mu)
 \leq24\eps^{-2}\norm{\mu}_{\WP}^2+
 6\sum_{\alpha\in\sigma_\eps(X)}
 \frac1{\ell_\alpha}\int_X|\mu|^2K_{\alpha,2}\dA \\
 &\leq 24\eps^{-2}\norm{\mu}_{\WP}^2+
 \sum_{\alpha\in\sigma_\eps(X)}
 \left(\frac{9\pi}{4}|c_{\alpha,0}|^2+6\frac{Q_\alpha^\perp(\mu)}{\ell_\alpha} \right).
\end{aligned}
\]
Here we apply \eqref{eq:unfold-Q} in the last inequality. By \eqref{eq:nonprincipal-inner} and
\eqref{eq:nonprincipal-outer}, the nonprincipal part satisfies
\[
 \sum_{\alpha\in\sigma_\eps(X)}
 \frac{Q_\alpha^\perp(\mu)}{\ell_\alpha}
 \leq(16\pi+76)\norm{\mu}_{\WP}^2
 <127\norm{\mu}_{\WP}^2.
\]
By \eqref{eq:firstvariation},
$|c_{\alpha,0}|^2
=4(a_\alpha^2+b_\alpha^2)/r_\alpha^2$. We consequently obtain
\[
 -\Ric_{\WP}(v,v)
 \leq(24\eps^{-2}+762)\norm{v}_{\WP}^2+
 9\pi\sum_{\alpha\in\sigma_\eps(X)}
 \frac{a_\alpha(v)^2+b_\alpha(v)^2}{r_\alpha^2}.
\]
Since $9\pi<30$ and $\eps^{-2}\geq100$, this is bounded by the
right-hand side of \eqref{eq:directional-Ricci}.
\end{proof}
The first term in \eqref{eq:directional-Ricci} is harmless once
\(\eps\) is fixed.  The second records the full singular behavior and
is precisely the term to be integrated along a minimizing geodesic.

For comparison, Wolpert's metric estimate
\cite[Theorem~4.3]{W08} gives, on every bounded-length Bers region
associated with a pants decomposition \(P\), $\|v\|_{\WP}^2$ is comparable to
 $ \sum_{\alpha\in P}
 \bigl(a_\alpha(v)^2+b_\alpha(v)^2\bigr)$, up to genus-dependent constant multiplication.
Thus Proposition~\ref{prop:Ricci} identifies \(r_\alpha^{-2}\) as the
additional curvature singularity in each short root-length direction.
The genus-uniform constants used below come from
Proposition~\ref{prop:Ricci}, rather than from this metric comparison.

\section{Twist deformations and the energy comparison}\label{sec:twists}
Let $\gamma:[0,L]\to\T_g$ be a unit-speed Weil--Petersson geodesic
whose projection to $\M_g$ minimizes distance, and suppose that
\[
 \min\left\{\ell_{\mathrm{sys}}(\gamma(0)),\ \ell_{\mathrm{sys}}(\gamma(L))\right\}\geq\eta>0.
\]
The argument has three parts. We first control the distortion of a
simultaneous normalized twist. We then twist only while a curve is
short, choosing an integral twist at the exit time. Finally, minimality
of $\pi(\gamma)$ turns the resulting comparison path into an estimate for
the angular energy.

\subsection{The normalized twist flow}
We retain the notation established earlier. Thus
$r_\alpha=\sqrt{\ell_\alpha}$,
$\lambda_\alpha=\grad r_\alpha$, and
\[
 a_\alpha(v)=\ip{\lambda_\alpha}{v},
 \qquad b_\alpha(v)=\ip{J\lambda_\alpha}{v}.
\]
The Fenchel--Nielsen twist vector is
$t_\alpha=\frac12J\grad\ell_\alpha$. Define
\begin{equation}\label{eq:Vdef}
 V_\alpha=r_\alpha^3J\lambda_\alpha
 =\ell_\alpha t_\alpha.
\end{equation}
For $s\in\R$, let
$\Phi_\alpha^s=\operatorname{Fl}_{V_\alpha}^s$ denote the time-$s$
flow of $V_\alpha$ on $\T_g$.
By Wolpert's twist--length duality and normalized angle coordinates
\cite[Section~4.2, before Corollary~4.8]{W08}, this flow is defined for
all $s\in\R$, preserves $\ell_\alpha$, and $\Phi_\alpha^1$ is the Dehn
twist about $\alpha$. Our $V_\alpha$ is $2\pi$ times the angle
variation denoted by $T_\alpha$ in that reference.

If $\sigma=\{\alpha_1,\ldots,\alpha_m\}$ is a multicurve and
$\theta\in\R^\sigma$, define
\begin{equation}\label{eq:Ftheta}
 F_\theta^\sigma(X)=
 \bigl(\Phi_{\alpha_1}^{\theta_{\alpha_1}}\circ\cdots\circ
 \Phi_{\alpha_m}^{\theta_{\alpha_m}}\bigr)(X).
\end{equation}
The factors in \eqref{eq:Ftheta} are diffeomorphisms of $\T_g$.
Because twist flows about disjoint curves commute and preserve all
the corresponding length functions, $F_\theta^\sigma$ is independent
of the enumeration of $\sigma$, and
\begin{equation}\label{eq:Faction}
 F_{\theta+\vartheta}^\sigma
 =F_\theta^\sigma\circ F_\vartheta^\sigma,\qquad
 F_0^\sigma=\mathrm{id}.
\end{equation}
For $k\in\Z^\sigma$, $F_k^\sigma$ is a product of Dehn twists and hence
a Weil--Petersson isometry.
\begin{lemma}\label{lem:infinitesimal}
Suppose that $r_\alpha\leq u\leq\sqrt{\ell_0}$ for all
$\alpha\in\sigma$. If $|\theta_\alpha|\leq1/2$ and
$W_\theta=\sum\limits_{\alpha\in\sigma}\theta_\alpha V_\alpha$, then
\begin{equation}\label{eq:infinitesimal}
 \bigl|(\mathcal L_{W_\theta}\gwp)(v,v)\bigr|
 \leq 81u^4\norm{v}_{\WP}^2.
\end{equation}
\end{lemma}

\begin{proof}
Since $J$ is parallel and $\lambda_\alpha=\grad r_\alpha$,
\[
 \nabla_vV_\alpha
 =3r_\alpha^2a_\alpha(v)J\lambda_\alpha
  +r_\alpha^3J\nabla_v\lambda_\alpha.
\]
Taking the inner product with $v$ and inserting the definition of
$H_\alpha$ gives the exact cancellation
\begin{equation}\label{eq:Lie}
 \tfrac12(\mathcal L_{V_\alpha}\gwp)(v,v)
 =3r_\alpha^2a_\alpha(v)b_\alpha(v)
  -r_\alpha^3\Hess r_\alpha(v,Jv)
 =-r_\alpha^3H_\alpha(v,Jv).
\end{equation}
Here we used $b_\alpha(Jv)=a_\alpha(v)$. By linearity of the Lie derivative and \eqref{eq:Lie},
\[
\begin{aligned}
 \bigl|(\mathcal L_{W_\theta}\gwp)(v,v)\bigr|
 =
 2\left|
 \sum_{\alpha\in\sigma}
 \theta_\alpha r_\alpha^3H_\alpha(v,Jv)
 \right|                                                                                
 \leq
 \sum_{\alpha\in\sigma}
 r_\alpha^3|H_\alpha(v,Jv)|,
\end{aligned}
\]
where the last inequality uses $|\theta_\alpha|\leq1/2$.
Since $H_\alpha\geq0$,
Cauchy--Schwarz for each positive semidefinite form and then for the
finite sum, together with $r_\alpha^3\leq u^4/r_\alpha$, yields
\begin{align*}
 \bigl|(\mathcal L_{W_\theta}\gwp)(v,v)\bigr|
 &\leq u^4
 \left(\sum_\alpha\frac{H_\alpha(v,v)}{r_\alpha}\right)^{1/2}
 \left(\sum_\alpha\frac{H_\alpha(Jv,Jv)}{r_\alpha}\right)^{1/2}\\
 &\leq 81u^4\norm{v}_{\WP}^2
\end{align*}
by \eqref{eq:Hsum}.
\end{proof}

\begin{proposition}
\label{prop:distortion}
Let $\sigma$ be a multicurve and let
$0<u\leq\sqrt{\ell_0}$. On the region where
$r_\alpha\leq u$ for every $\alpha\in\sigma$, one has
\begin{equation}\label{eq:distortion}
   e^{-81u^4}\gwp
   \leq
   (F_\theta^\sigma)^*\gwp
   \leq
   e^{81u^4}\gwp
   \qquad
   (\theta\in\R^\sigma).
\end{equation}
\end{proposition}

\begin{proof}
First suppose that $|\theta_\alpha|\leq1/2$ for every
$\alpha\in\sigma$. Since the vector fields $V_\alpha$ commute,
$\Psi_s=F_{s\theta}^{\sigma}$ is the flow of
$W_\theta=\sum_\alpha\theta_\alpha V_\alpha$. For $X\in\T_g$ and
$0\ne v\in T_X\T_g$, set
\[
 q(s)=(\Psi_s^*\gwp)_X(v,v).
\]
Every $r_\alpha$ is preserved by $\Psi_s$. Lemma~\ref{lem:infinitesimal},
applied at $\Psi_s(X)$ to $(\mathrm d\Psi_s)_Xv$, therefore gives
\[
 |q'(s)|
 =\left|\bigl(\mathcal L_{W_\theta}\gwp\bigr)_{\Psi_s(X)}
   \bigl((\mathrm d\Psi_s)_Xv,(\mathrm d\Psi_s)_Xv\bigr)\right|
 \leq81u^4q(s).
\]
Integrating $| (\log q)'|\leq81u^4$ from $0$ to $1$ proves
\eqref{eq:distortion} when $|\theta_\alpha|\leq1/2$.

For arbitrary $\theta\in\R^\sigma$, choose $k\in\Z^\sigma$ such that
\[
   \theta_0=\theta-k\in[-1/2,1/2]^\sigma.
\]
The group law gives
$F_\theta^\sigma=F_k^\sigma\circ F_{\theta_0}^\sigma$. Since
$F_k^\sigma$ is a Weil--Petersson isometry,
\[
 (F_\theta^\sigma)^*\gwp=(F_{\theta_0}^\sigma)^*\gwp,
\]
and the preceding estimate applies to $\theta_0$.
\end{proof}

In Proposition~\ref{prop:distortion}, the parameters are fixed while a
map is differentiated. Later the parameters will depend on time. Their
integer reductions will only be used to invoke this proposition; they
will not be differentiated.

\subsection{Comparison paths} The purpose of this subsection is to use the normalized twist maps from the preceding subsection to construct competitors for \(\pi(\gamma)\).  On each short-curve
excursion \(I_\alpha(u)\), the twist parameter is chosen to vanish at
entry and to be integral at exit; hence the resulting path is
continuous in \(\M_g\) and has the same endpoints as \(\pi(\gamma)\). In this subsection, we compute the velocity of this comparison path. By choosing the twisting functions to
cancel the angular part of \(\dot\gamma\), quotient minimality then
yields the angular-energy estimate of the next subsection.

Recall that $\gamma:[0,L]\to\T_g$ is a unit-speed Weil--Petersson geodesic
whose projection minimizes distance in $\M_g$, with both endpoints of
systole at least $\eta>0$. We assume $L>0$, since the case $L=0$ is
immediate. Along $\gamma$, write
\[
 r_\alpha(t)=r_\alpha(\gamma(t)),\quad
 a_\alpha(t)=a_\alpha(\dot\gamma(t))=r_\alpha'(t),\quad
 b_\alpha(t)=b_\alpha(\dot\gamma(t)).
\]
For $0<u<\sqrt\eta$, put
\begin{equation}\label{eq:Idef}
 I_\alpha(u)=\{t\in[0,L]:r_\alpha(t)\leq u\}.
\end{equation}

\noindent Recall that the root length $r_\alpha$ is convex along Weil--Petersson geodesics
\cite[Corollary~3.4 and Example~3.5]{W08}. Thus each nondegenerate
visit is an interval
$I_\alpha(u)=[p_\alpha(u),q_\alpha(u)]\subset(0,L)$,
with
\begin{equation}\label{eq:interval-endpoints}
 r_\alpha(p_\alpha(u))=r_\alpha(q_\alpha(u))=u,\qquad
 a_\alpha(p_\alpha(u))\leq0\leq a_\alpha(q_\alpha(u)).
\end{equation}
Only finitely many curves have nonempty excursions. Indeed,
compactness of $\gamma([0,L])$ gives a uniform bilipschitz comparison
with $\gamma(0)$, and hence a constant $C_\gamma$ such that
$\ell_\alpha(\gamma(0))\leq C_\gamma\ell_\alpha(\gamma(t))$ for every
$\alpha$ and $t$. A curve with $I_\alpha(u)\ne\varnothing$ therefore
has length at most $C_\gamma u^2$ on $\gamma(0)$, and a fixed surface
has only finitely many closed geodesics of bounded length.
The constant $C_\gamma$ is used only for this finiteness assertion.

Assume also $u\leq\sqrt{\ell_0}$. The curves active at a given time
are disjoint by the collar lemma. Let $\mathcal A$ be any subfamily of
the nondegenerate visits, and choose smooth functions $h_\alpha$ on
$I_\alpha(u)$ such that
\begin{equation}\label{eq:integer-exits}
 h_\alpha(p_\alpha(u))=0,\qquad
 h_\alpha(q_\alpha(u))\in\Z.
\end{equation}
For a time away from the finitely many interval endpoints, define
\[
 \sigma_{\mathcal A}(t)
 =\{\alpha\in\mathcal A:t\in I_\alpha(u)\},
 \qquad
 h(t)=(h_\alpha(t))_{\alpha\in\sigma_{\mathcal A}(t)}.
\]
Set
\begin{equation}\label{eq:competitor}
 \zeta(t)=\pi\bigl(
 F_{h(t)}^{\sigma_{\mathcal A}(t)}(\gamma(t))
 \bigr).
\end{equation}
For an empty active multicurve, this means
$\zeta(t)=\pi(\gamma(t))$.
\begin{lemma}
\label{lem:comparison-path}
The curve $\zeta$ defined in \eqref{eq:competitor} extends uniquely to a
continuous, piecewise smooth path in $\M_g$ with the same endpoints as
$\pi(\gamma)$. On every open interval on which
$\sigma_{\mathcal A}(t)=\sigma$ is constant, it has the local lift
\[
   \widetilde\zeta(t)
   =
   F_{h(t)}^\sigma(\gamma(t)),
\]
whose velocity is
\begin{equation}\label{eq:velocity}
   \dot{\widetilde\zeta}(t)
   =
   (\mathrm dF_{h(t)}^\sigma)_{\gamma(t)}
   \left(
      \dot\gamma(t)
      +
      \sum_{\alpha\in\sigma}
      h_\alpha'(t)V_\alpha(\gamma(t))
   \right).
\end{equation}
\end{lemma}

\begin{proof}
At an entry time the new twist factor is the identity. At an exit
time it is an integral Dehn twist, so it disappears after projection.
All curves active at a transition have length at most $\ell_0$ and
are disjoint; their twist flows commute. Thus the one-sided limits
agree even at simultaneous entries and exits. There are only finitely
many transitions, and none occurs at $0$ or $L$ because
$u<\sqrt\eta$. This proves continuity, piecewise smoothness, and
agreement of the endpoints.

On an interval with fixed active multicurve, commutativity gives
$(\mathrm dF_h^\sigma)V_\alpha=V_\alpha\circ F_h^\sigma$.
The chain rule now gives
\[
\begin{aligned}
   \dot{\widetilde\zeta}(t)
   &=
   (\mathrm dF_{h(t)}^\sigma)_{\gamma(t)}
   \dot\gamma(t)
   +
   \sum_{\alpha\in\sigma}
   h_\alpha'(t)
   V_\alpha\bigl(\widetilde\zeta(t)\bigr)                 \\
   &=
   (\mathrm dF_{h(t)}^\sigma)_{\gamma(t)}
   \left(
      \dot\gamma(t)
      +
      \sum_{\alpha\in\sigma}
      h_\alpha'(t)V_\alpha(\gamma(t))
   \right),
\end{aligned}
\]
which proves \eqref{eq:velocity}.
\end{proof}

\subsection{The energy inequality}
For a nondegenerate excursion define
\begin{equation}\label{eq:Mdef}
 M_\alpha(u)=\int_{I_\alpha(u)}r_\alpha^{-6}\dd.
\end{equation}
This number is positive and finite because $r_\alpha>0$ on the compact
interval $I_\alpha(u)$.
\begin{lemma}\label{lem:integer-correction}
For every $\alpha\in\mathcal A$, there are smooth functions
$h_\alpha,e_\alpha:I_\alpha(u)\to\R$, with $h_\alpha$ satisfying
\eqref{eq:integer-exits}, such that, for almost every $t$,
\begin{align}
 \norm{\dot\gamma+\sum_{\alpha\in\sigma_{\mathcal A}(t)}
 h_\alpha'V_\alpha}_{\WP}^2
 &\leq1-\sum_{\alpha\in\sigma_{\mathcal A}(t)}b_\alpha^2
 +\sum_{\alpha\in\sigma_{\mathcal A}(t)}e_\alpha^2,
 \label{eq:square-completion}\\
 \int_{I_\alpha(u)}e_\alpha^2\dd
 &\leq\frac1{4M_\alpha(u)}.\label{eq:rounding-cost}
\end{align}
\end{lemma}

\begin{proof}
Choose an integer $k_\alpha$ nearest to
$-\int_{I_\alpha(u)}b_\alpha r_\alpha^{-3}\dd$, and define
\[
 \kappa_\alpha M_\alpha(u)
 =k_\alpha+\int_{I_\alpha(u)}b_\alpha r_\alpha^{-3}\dd.
\]
Thus $|\kappa_\alpha M_\alpha(u)|\leq1/2$. Set
\begin{equation}\label{eq:twist-correction}
 h_\alpha(t)=\int_{p_\alpha(u)}^t
 \bigl(-b_\alpha(s)r_\alpha(s)^{-3}
       +\kappa_\alpha r_\alpha(s)^{-6}\bigr)\,ds,
 \qquad e_\alpha=\kappa_\alpha r_\alpha^{-3}.
\end{equation}
Then $h_\alpha(p_\alpha(u))=0$,
$h_\alpha(q_\alpha(u))=k_\alpha \in \Z$, and
\begin{equation}\label{eq:hprime}
 h_\alpha'=-b_\alpha r_\alpha^{-3}+\kappa_\alpha r_\alpha^{-6},
 \qquad h_\alpha'V_\alpha=(-b_\alpha+e_\alpha)J\lambda_\alpha.
\end{equation}
At a time away from the excursion endpoints, all sums below are over
$\sigma_{\mathcal A}(t)$. By \eqref{eq:dual-Bessel},
\[
 \begin{aligned}
 \norm{\dot\gamma+\sum_\alpha h_\alpha'V_\alpha}_{\WP}^2
 &= 1+2\sum_\alpha(-b_\alpha+e_\alpha)b_\alpha
          +\left\|\sum_\alpha (-b_\alpha+e_\alpha) J\lambda_\alpha
\right\|_{\mathrm{WP}}^2\\
 &\leq1+2\sum_\alpha(-b_\alpha+e_\alpha)b_\alpha
          +\sum_\alpha(-b_\alpha+e_\alpha)^2\\
 &=1-\sum_\alpha b_\alpha^2+\sum_\alpha e_\alpha^2.
 \end{aligned}
\]
Finally, we have
$$\int_{I_\alpha(u)}e_\alpha^2\dd
=\kappa_\alpha^2M_\alpha(u)\leq1/(4M_\alpha(u))$$
as desired.
\end{proof}

\begin{proposition}\label{prop:energy}
If $0<u\leq\sqrt{\ell_0}$ and $u<\sqrt\eta$, then
\begin{equation}\label{eq:energy}
   \sum_{\alpha\in\mathcal A}
   \int_{I_\alpha(u)} b_\alpha(t)^2\,dt
   \leq
   \frac14
   \sum_{\alpha\in\mathcal A}M_\alpha(u)^{-1}
   +
   81u^4L.
\end{equation}
\end{proposition}

\begin{proof}
Use the functions from Lemma~\ref{lem:integer-correction} in
\eqref{eq:competitor}, and set
$E(\zeta)=\int_0^L|\dot\zeta|^2\,dt$. Proposition~\ref{prop:distortion}
and \eqref{eq:square-completion} give
\[
 E(\zeta)\leq e^{81u^4}\left(
 L-\sum_{\alpha\in\mathcal A}\int_{I_\alpha(u)}b_\alpha^2\dd
 +\frac14\sum_{\alpha\in\mathcal A}M_\alpha(u)^{-1}\right).
\]
The projected geodesic $\pi(\gamma)$ minimizes between its endpoints and
has length $L$. Since $\zeta$ has the same endpoints,
Cauchy--Schwarz gives
\[
 E(\zeta)\geq\frac{\operatorname{Length}(\zeta)^2}{L}\geq L.
\]
Combining the two displays yields \eqref{eq:energy}, because
$1-e^{-81u^4}\leq 81u^4$.
\end{proof}

\begin{remark}
The weight in \eqref{eq:hprime} is forced by the least-cost correction.
For any prescribed value of $\int_Iq\,dt$,
\[
 \left(\int_Iq\,dt\right)^2
 \leq\left(\int_Ir^6q^2\,dt\right)
      \left(\int_Ir^{-6}\,dt\right),
\]
with equality precisely when $q$ is proportional to $r^{-6}$.
Thus rounding by a discrepancy $\delta_\alpha\in[-1/2,1/2]$
has minimum cost $\delta_\alpha^2/M_\alpha(u)$, at most
$1/(4M_\alpha(u))$.
\end{remark}

\section{Integrating the short-curve contribution}\label{sec:integrals}

We retain the minimizing segment and the excursion intervals introduced
in Section~\ref{sec:twists}. Set
\begin{equation}\label{eq:rhochoice}
 \rho=\frac14\min\{\sqrt{\ell_0},\sqrt\eta\}.
\end{equation}
The factor $1/4$ leaves room to compare the threshold $s$ with the
larger threshold $2s$. At both scales the excursions remain away from
the endpoints of $\gamma$, and the active curves form a multicurve.

\subsection{An angular estimate at one scale}

For $0<s\leq\rho$, set
\[
 \mathcal E(s)=\sum_\alpha\int_{I_\alpha(s)}b_\alpha^2\dd.
\]
All sums over visits are finite, as observed in
Section~\ref{sec:twists}.

\begin{figure}[htbp]
\centering
\begin{tikzpicture}[x=0.90cm,y=0.90cm,>=Latex,font=\small]
 \definecolor{curveblue}{RGB}{39,70,102}
 \definecolor{bandgold}{RGB}{165,111,45}
 \fill[blue!4] (1.382,0.55) rectangle (7.018,2.55);
 \fill[orange!10] (2.351,0.55) rectangle (6.049,1.55);
 \draw[->] (0.45,0.55)--(8.0,0.55) node[right] {$t$};
 \draw[->] (0.65,0.15)--(0.65,3.8) node[above] {$r_\alpha(t)$};
 \draw[dashed,gray] (0.65,1.55)--(7.55,1.55);
 \draw[dashed,gray] (0.65,2.55)--(7.55,2.55);
 \node[left] at (0.65,1.55) {$s$};
 \node[left] at (0.65,2.55) {$2s$};
 \draw[very thick,curveblue,domain=0.95:7.45,samples=100]
 plot (\x,{0.8+0.22*(\x-4.2)*(\x-4.2)});
 \draw[densely dotted] (1.382,0.55)--(1.382,2.55);
 \draw[densely dotted] (7.018,0.55)--(7.018,2.55);
 \fill[curveblue] (1.382,2.55) circle (1.5pt);
 \fill[curveblue] (7.018,2.55) circle (1.5pt);
 \draw[thin,gray] (1.382,2.55)--(2.35,3.18)
 node[above,text=black] {$h_\alpha=0$};
 \draw[thin,gray] (7.018,2.55)--(6.05,3.18)
 node[above,text=black] {$h_\alpha\in\Z$};
 \draw[<->,bandgold,thick] (2.351,-0.05)--(6.049,-0.05)
 node[midway,fill=white,inner sep=3pt] {$I_\alpha(s)$};
 \draw[<->,curveblue,thick] (1.382,-0.6)--(7.018,-0.6)
 node[midway,fill=white,inner sep=3pt] {$I_\alpha(2s)$};
\end{tikzpicture}
\caption{The twist comparison
is made on the outer interval. Every visit to the inner interval forces
radial variation at least $2s$ on the outer one. The integer twist at exit
preserves the endpoint in moduli space.}
\label{fig:threshold}
\end{figure}
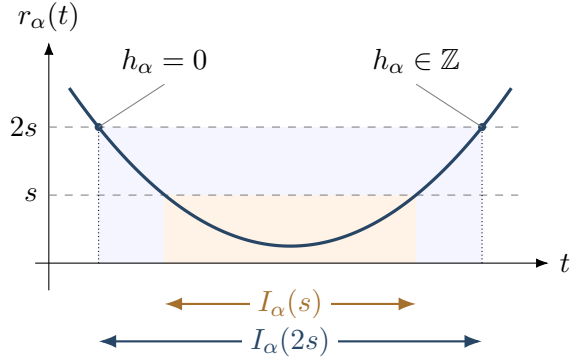

\begin{lemma}\label{lem:angular-scale}
For every $0<s\leq\rho$,
\begin{equation}\label{eq:angular-scale}
 \mathcal E(s)\leq 1300s^4L.
\end{equation}
\end{lemma}

\begin{proof}
Let $\mathcal A_s$ be the curves with $\min r_\alpha\leq s$.
For $\alpha\in\mathcal A_s$, write
$I_\alpha(2s)=[p_\alpha,q_\alpha]$,
$T_\alpha=q_\alpha-p_\alpha$, and $m_\alpha=\min r_\alpha$.
So $r_\alpha(p_\alpha)
=r_\alpha(q_\alpha)=2s$. Let $[c_\alpha^-,c_\alpha^+]$ be the set on which $r_\alpha$
attains its minimum. Convexity of $r_\alpha$ implies that
$r_\alpha'\leq0$ on $[p_\alpha,c_\alpha^-]$,
$r_\alpha'=0$ on $[c_\alpha^-,c_\alpha^+]$, and
$r_\alpha'\geq0$ on $[c_\alpha^+,q_\alpha]$. Since
$a_\alpha=r_\alpha'$, it follows that
\[
\begin{aligned}
   \int_{I_\alpha(2s)}|a_\alpha(t)|\,dt
   &=
   r_\alpha(p_\alpha)-m_\alpha
   +
   r_\alpha(q_\alpha)-m_\alpha  \\
   &=
   2(2s-m_\alpha)
   \geq2s.
\end{aligned}
\]
Cauchy--Schwarz and $r_\alpha\leq2s$ on this interval imply
\[
 \int_{I_\alpha(2s)}a_\alpha^2\dd\geq\frac{4s^2}{T_\alpha},
 \qquad
 M_\alpha(2s)\geq\frac{T_\alpha}{(2s)^6}.
\]
Consequently,
\begin{equation}\label{eq:rounding-scale}
 M_\alpha(2s)^{-1}
 \leq16s^4\int_{I_\alpha(2s)}a_\alpha^2\dd.
\end{equation}
The curves with $t\in I_\alpha(2s)$ form a multicurve.
Summing \eqref{eq:rounding-scale} over $\alpha\in\mathcal A_s$ and applying \eqref{eq:Bessel} pointwise yields
\[
\begin{aligned}
\sum_{\alpha\in\mathcal A_s}M_\alpha(2s)^{-1}
&\le 16s^4\sum_{\alpha\in\mathcal A_s}
\int_0^L
\mathbf{1}_{I_\alpha(2s)}(t)\,a_\alpha(t)^2\,dt \\
&=16s^4\int_0^L
\sum_{\substack{\alpha\in\mathcal A_s\\
t\in I_\alpha(2s)}}
a_\alpha(t)^2\,dt\\
&\leq16s^4L.
\end{aligned}
\]
Apply Proposition~\ref{prop:energy} to $\mathcal A_s$ at threshold
$2s$. Its hypotheses hold by \eqref{eq:rhochoice}, and
$I_\alpha(s)\subset I_\alpha(2s)$. Therefore we get
\[
 \mathcal E(s)\leq\tfrac14(16s^4L)+81(2s)^4L=1300s^4L
\]
as desired.
\end{proof}

The larger threshold guarantees radial variation at least $2s$, even
when an excursion barely reaches $s$. This is what controls the
rounding cost; see Figure~\ref{fig:threshold}.

\subsection{Proof of the integrated Ricci bound}

We now prove Theorem~\ref{thm:integrated-intro}. Define the weighted
angular and radial energies by
\begin{equation}\label{eq:ABdef}
   B_\rho
   =
   \sum_\alpha
   \int_{I_\alpha(\rho)}
   \frac{b_\alpha(t)^2}{r_\alpha(t)^2}\,dt,
   \qquad
   A_\rho
   =
   \sum_\alpha
   \int_{I_\alpha(\rho)}
   \frac{a_\alpha(t)^2}{r_\alpha(t)^2}\,dt.
\end{equation}

\begin{lemma}
\label{lem:weighted-energy}
The angular and radial contributions satisfy
\begin{align}
   B_\rho
   &\leq
   2600\rho^2L,
   \label{eq:angular}\\
   A_\rho
   &\leq
   3B_\rho+81L,
   \label{eq:radial-one}\\
   A_\rho+B_\rho
   &\leq
   \bigl(81+10400\rho^2\bigr)L.
   \label{eq:radial}
\end{align}
\end{lemma}

\begin{proof}
For $0<r\leq\rho$,
$r^{-2}=\rho^{-2}+2\int_r^\rho s^{-3}\,ds$.
Applying this identity with
$r=r_\alpha(t)$, multiplying by $b_\alpha(t)^2$, and integrating over
$I_\alpha(\rho)$ gives
\begin{equation}\label{eq:layer-cake-energy}
 B_\rho=\rho^{-2}\mathcal E(\rho)
       +2\int_0^\rho s^{-3}\mathcal E(s)\,ds.
\end{equation}
Lemma~\ref{lem:angular-scale} therefore yields
\[
 B_\rho\leq1300\rho^2L+2600L\int_0^\rho s\,ds
 =2600\rho^2L.
\]

For the radial term,  since
$a_\alpha=r_\alpha'$ and $\gamma$ is a geodesic,
$ a_\alpha'
   =
   \Hess r_\alpha(\dot\gamma,\dot\gamma).$
By the definition of $H_\alpha$, $
   \Hess r_\alpha(\dot\gamma,\dot\gamma)
   =
   \frac{3b_\alpha^2}{r_\alpha}
   +
   H_\alpha(\dot\gamma,\dot\gamma).$
Therefore
\begin{equation}\label{eq:radial-identity}
 \left(\frac{a_\alpha}{r_\alpha}\right)'
 =\frac{3b_\alpha^2-a_\alpha^2}{r_\alpha^2}
  +\frac{H_\alpha(\dot\gamma,\dot\gamma)}{r_\alpha}.
\end{equation}
On a nondegenerate interval $I_\alpha(\rho)=[p_\alpha,q_\alpha]$,
\eqref{eq:interval-endpoints} gives
\[
 \left[\frac{a_\alpha}{r_\alpha}\right]_{p_\alpha}^{q_\alpha}
 =\frac{a_\alpha(q_\alpha)-a_\alpha(p_\alpha)}{\rho}\geq0.
\]
Integrating \eqref{eq:radial-identity}, summing over the excursions,
and applying \eqref{eq:Hsum} to the active multicurve gives
$$A_\rho\leq3B_\rho+81L.$$ Empty or degenerate intervals contribute
nothing. The final bound follows by adding $B_\rho$ and using
\eqref{eq:angular}.
\end{proof}

\begin{proof}[Proof of Theorem~\ref{thm:integrated-intro}]
Apply Proposition~\ref{prop:Ricci} with $\eps=\rho^2$.
Since $\ell_\alpha\leq\rho^2$ is equivalent to
$t\in I_\alpha(\rho)$, integration gives
\begin{equation}\label{eq:Ricci-AB}
 \int_0^L-\Ric_{\WP}(\dot\gamma,\dot\gamma)\dd
 \leq32\rho^{-4}L+30(A_\rho+B_\rho).
\end{equation}
Lemma~\ref{lem:weighted-energy} now gives
\begin{equation}\label{eq:numerical-integrated}
\begin{aligned}
   \int_0^L
   -\Ric_{\WP}(\dot\gamma,\dot\gamma)\,dt
   &\leq
   \left[
      32\rho^{-4}
      +
      30\bigl(81+10\,400\rho^2\bigr)
   \right]L                                                \\
   &=
   \bigl(
      32\rho^{-4}
      +2430
      +312\,000\rho^2
   \bigr)L\\
&\leq
   C(1+\eta^{-2})L.
\end{aligned}
\end{equation}
Using \eqref{eq:rhochoice}, this proves the result with suitable \(C\geq 1\). More precisely, we may choose
\[
 A_\eta=32\rho^{-4}+2430+312\,000\rho^2.
\]
In particular, when $\eta=1$, we have $\rho^2=1/160$ and $A_1=819\,200+2430+1950=823\,580.$
\end{proof}


\section{Thick endpoints of metric balls}\label{sec:balls}

We translate Theorem~\ref{thm:integrated-intro} into a volume estimate under a standard argument. The first step is a radial Jacobian inequality which requires only an integrated Ricci bound. We then integrate this inequality over the part
of the polar domain which parametrizes minimizers with thick endpoints. 

\subsection{A radial Jacobian estimate}

Let $N$ be a Riemannian manifold of dimension $d\geq2$ with nonpositive
sectional curvature. Completeness is not assumed. Along a unit-speed
radial geodesic $\gamma:[0,t]\to N$, let $\mathcal J(s)$ denote the
polar volume density based at $\gamma(0)$. Thus
$\mathcal J(s)/s^{d-1}\to1$ as $s\downarrow0$. This density is positive
while the geodesic is defined, since there are no conjugate points.

\begin{lemma}\label{lem:Jacobi}
Set $m=d-1$ and $q(s)=-\Ric(\dot\gamma(s),\dot\gamma(s))$.
For every $\kappa>0$,
\begin{equation}\label{eq:Jacobi}
 \mathcal J(t)\leq\kappa^{-m}
 \exp\left(\frac{m\kappa t}{2}
       +\frac1{2\kappa}\int_0^tq(s)\,ds\right).
\end{equation}
\end{lemma}

\begin{proof}
Let $\mathbb{A}(s)$ be the radial Jacobi tensor in a parallel
orthonormal frame, with $\mathbb{A}(0)=0$ and $\mathbb{A}'(0)=I$.
For $s>0$, set
\[
 \mathcal J(s)=\det \mathbb{A}(s),\qquad
 S(s)= \mathbb{A}'(s)\ \mathbb{A}(s)^{-1},\qquad
 h=\operatorname{tr}S=\mathcal J'/\mathcal J.
\]
The Riccati equation and $\operatorname{tr}(S^2)\geq h^2/m$ imply
that $y=\mathcal J^{1/m}$ satisfies
\begin{equation}\label{eq:y-ineq}
 y''\leq\frac qm y,\qquad y(0)=0,\quad y'(0)=1;
\end{equation}
see \cite[Example~4.3]{Ball16}.
The initial data follow from the Jacobi equation:
$\mathbb{A}(s)=sI+O(s^3)$, so $y(s)=s+O(s^3)$, with the
corresponding differentiated expansion.

Nonpositive sectional curvature makes the radial shape operator
positive semidefinite. To see this directly, let $U$ be a radial Jacobi
field with $U(0)=0$. Integration by parts in the Jacobi equation gives
\[
 \langle U'(s),U(s)\rangle
 =\int_0^s\left(
 |U'|^2-\langle R(U,\dot\gamma)\dot\gamma,U\rangle\right)d\tau.
\]
The integrand is nonnegative when sectional curvature is nonpositive,
so $\langle S (U),U\rangle\geq0$. Hence $S$ is positive semidefinite,
$h\geq0$, and therefore $y'\geq0$.

Now fix $\kappa>0$ and set $F=(y')^2+\kappa^2y^2$. Since
$y,y'\geq0$, equation \eqref{eq:y-ineq} gives
\[
 \begin{split}
 F'\leq2\left(\frac qm+\kappa^2\right)yy'
 &=\left(\kappa+\frac{q}{m\kappa}\right)(2\kappa yy')\\
 &\leq\left(\kappa+\frac{q}{m\kappa}\right)F,
 \end{split}
\]
where the last step is $2\kappa yy'\leq(y')^2+\kappa^2y^2$.
Because $F(0)=1$, Gr\"onwall's inequality yields
\[
 F(t)\leq
 \exp\left(\kappa t+\frac1{m\kappa}\int_0^tq(s)\,ds\right).
\]
Finally $\kappa^2y(t)^2\leq F(t)$. Raising the resulting estimate for
$y(t)$ to the $m$-th power, and using $\mathcal J=y^m$, proves
\eqref{eq:Jacobi}. 
\end{proof}

\begin{remark}
Every step takes place along the given radial
geodesic, so no completeness assumption is required in Lemma \ref{lem:Jacobi}. If $\int_0^tq\leq At$ with $A>0$, choosing
$\kappa=\sqrt{A/m}$ gives
\begin{equation}\label{eq:Jacobi-optimized}
 \mathcal J(t)\leq(m/A)^{m/2}e^{\sqrt{Am}\,t}.
\end{equation}
\end{remark}

\subsection{The minimizing polar domain}

Fix $X\in\M_g^{\geq\eta}$ and a lift $\widetilde X\in\T_g$.
Let $S_{\widetilde X}$ be the unit sphere in $T_{\widetilde X}\T_g$.
For each $v\in S_{\widetilde X}$ let $\gamma_v$ be the radial geodesic
with initial velocity $v$, on its maximal interval of existence in
$\T_g$. For $R>0$, define
\begin{equation}\label{eq:Omega}
\begin{aligned}
 \Omega_R=\{(t,v):{}&0<t<R,\ \gamma_v(t)\text{ is defined},\\
 &\pi(\gamma_v(t))\in\M_g^{\geq\eta},\quad
 \dwp(X,\pi(\gamma_v(t)))=t\}.
\end{aligned}
\end{equation}
The set $\Omega_R$ is Borel. Indeed, the exponential map has open
domain, and endpoint thickness and the distance equality are closed
conditions in that domain. Every point of
$B_{\M_g}(X,R)\cap\M_g^{\geq\eta}$, other than $X$, has a preimage in
$\Omega_R$ by the interior minimizing-path property: choose a
distance-realizing segment in moduli space, lift it from
$\widetilde X$, and use nonrefraction to keep the lifted segment in
$\T_g$. Notice that no assertion is made that $\Omega_R$ is
star-shaped in the radial variable.

\begin{proof}[Proof of Theorem~\ref{thm:ball-intro}]
The real dimension is $d=6g-6$; put $m=6g-7$. For every
$(t,v)\in\Omega_R$, Theorem~\ref{thm:integrated-intro} and
the optimized bound \eqref{eq:Jacobi-optimized} imply
\begin{equation}\label{eq:Jthick}
 \mathcal J(t,v)\leq(m/A_\eta)^{m/2}e^{\sqrt{A_\eta m}\,t}.
\end{equation}
Indeed, both endpoints of $\gamma_v|_{[0,t]}$ project to the
$\eta$-thick part, and the defining distance equality in
\eqref{eq:Omega} says that its projection is globally minimizing.
Thus the hypotheses of Theorem~\ref{thm:integrated-intro} hold on this
particular segment.

The area formula for $\pi\circ\exp_{\widetilde X}$ gives
\begin{align*}
 \Vol\bigl(B_{\M_g}(X,R)\cap\M_g^{\geq\eta}\bigr)
 &\leq\int_{\Omega_R}\mathcal J(t,v)\,dt\,dS(v)\\
 &\leq |S^m|\,(m/A_\eta)^{m/2}
 \int_0^R e^{\sqrt{A_\eta m}\,t}\,dt.
\end{align*}
For the first inequality, take a countable cover of the open
exponential domain by relatively compact coordinate neighborhoods,
and partition $\Omega_R$ into disjoint Borel subsets subordinate to
this cover. Apply the area formula on each subset, using orbifold
charts on the target, and sum. Multiple preimages can only increase
the integral, while finite quotient multiplicities can only decrease
the image volume. The effective orbifold singular set has measure zero.

By Stirling's formula,
\[
 |S^m|\,m^{m/2}
 =\frac{2\pi^{(m+1)/2}m^{m/2}}{\Gamma((m+1)/2)}
 \leq e^{Cm}.
\]
Since $A_\eta\geq1$ and $m=6g-7\geq5$, the time integral is at most
$e^{\sqrt{A_\eta m}R}$: if $a=\sqrt{A_\eta m}>1$, then
$\int_0^Re^{at}\,dt=(e^{aR}-1)/a\leq e^{aR}$.
Also $A_\eta^{-m/2}\leq1$. Absorbing the sphere factor into
$\exp(Cg)$ and using $m<6g$ proves \eqref{eq:intro-ball}.
\end{proof}

Only the upper bound for the Jacobian is integrated over the full
product $(0,R)\times S_{\widetilde X}$. The integrated Ricci estimate
is used solely on $\Omega_R$, which need not be star-shaped.

\section{Proof of the diameter theorem}\label{sec:diameter}
We first provide an effective bound for the proportion of the Weil-Petersson volume of the $1$-thick part, which is essential in the proof of Theorem~\ref{thm:main}. One may also see \cite[Theorem 4.2]{Mir13} for implicit lower and upper bounds.
\begin{proposition}\label{prop:thickvolume}
For all sufficiently large $g$,
\begin{equation}\label{eq:thickvolume}
 \Vol(\M_g^{\geq1})\geq\tfrac13\Vol(\M_g)
 \geq\exp(2g\log g-Cg),
\end{equation}
where $C$ is universal.
\end{proposition}
\begin{proof}
Let $\widehat V_{h,k}(\mathbf L)$ denote the volumes in
\cite{Mir13}, retaining its boundary-labeling and orbifold
conventions, including the one-handle factors in the integration
formula. Its symplectic form satisfies
$\widehat\omega_{\mathrm{WP}}=2\omega_{\mathrm{WP}}.$
Hence, for closed surfaces of sufficiently large genus,
\[
d\widehat V_g
=2^{3g-3}\,d\operatorname{Vol}_{\mathrm{WP}},
\qquad
\operatorname{Vol}_{\mathrm{WP}}(\mathcal M_g)
=2^{-(3g-3)}\widehat V_g.
\]
The normalized probability measures therefore agree.
All cutting formulas and cross-dimensional volume ratios below
are used in the hatted normalization.

Write $\mathbb P_g$
for normalized WP measure and $\mathbb E_g$ for its expectation.
Let $N(X)$ count simple closed geodesics of length less than $1$.
Since $\{\ell_{\mathrm{sys}}<1\}=\{N\geq1\}$, Markov's inequality gives
\[
 \mathbb P_g(\ell_{\mathrm{sys}}<1)\leq\mathbb E_g[N].
\]

The integration formula \cite[Theorem~2.2]{Mir13}, with symmetry
factors bounded by one, gives
\[
 \mathbb P_g(\ell_{\mathrm{sys}}<1)\leq\mathbb E_g[N]
 \leq\frac1{\widehat V_g}\int_0^1 t
 \left(\widehat V_{g-1,2}(t,t)
 +\sum_{i=1}^{\lfloor g/2\rfloor}
 \widehat V_{i,1}(t)\widehat V_{g-i,1}(t)\right)\,dt.
\]
The first term counts nonseparating curves and the sum counts
separating curves. 

The coefficient bound and volume polynomial in
\cite[Lemma~3.2(2) and (3.1)]{Mir13} imply
\[
 \widehat V_{h,k}(\mathbf L)
 \leq\widehat V_{h,k}\prod_{j=1}^k\cosh(L_j/2).
\]
Moreover, \cite[Theorem~1.2 and (3.19)]{Mir13} gives
\[
 \frac{\widehat V_{g-1,2}}{\widehat V_g}=1+O(g^{-1}),
 \qquad
 \frac{\sum_{i=1}^{\lfloor g/2\rfloor}
 \widehat V_{i,1}\widehat V_{g-i,1}}{\widehat V_g}=O(g^{-1}).
\]
For $0\leq t\leq1$, the nonseparating contribution divided by
$\widehat V_g$ is consequently at most
$(1+O(g^{-1}))\cosh^2(t/2)$, uniformly in $t$, and the total
separating contribution is $O(g^{-1})\cosh^2(t/2)$. Hence
\[
 \mathbb P_g(\ell_{\mathrm{sys}}<1)
 \leq\int_0^1t\cosh^2(t/2)\,dt+O(g^{-1})
 =\frac34-\frac1{2e}+O(g^{-1})<\frac23
\]
for large $g$. This proves the first inequality in \eqref{eq:thickvolume}.

Finally,
\cite[Corollary~3.6]{Mir13} and Stirling's formula give
$\log\widehat V_g=2g\log g+O(g)$. Passing to our volume normalization
subtracts $(3g-3)\log2$, and the factor $1/3$ subtracts only $\log3$.
Both changes are absorbed by the $O(g)$ term, proving the second
inequality.
\end{proof}

The proof of Theorem~\ref{thm:main} is now a direct comparison of the large thick-part volume with the growth of thick endpoints in a metric ball. Write $$D_{g,n}=\diam(\M_{g,n}) \quad \textit{$\mathrm{and}$} \quad D_g=\diam(\M_{g}).$$ Passing to the Weil--Petersson
completion, denoted by a bar, does not change the diameter. We first recall the following property that transfers the estimate for closed surfaces to punctured surfaces. 
\begin{lemma}\label{lem:puncture-transfer}
For every fixed $n\geq1$,
\[
 \liminf_{g\to\infty}\frac{D_{g,n}}{\sqrt g\log g}
 \geq
 \liminf_{g\to\infty}\frac{D_{g}}{\sqrt g\log g}.
\]
\end{lemma}

\begin{proof}
This is essentially \cite[Lemma~6.1]{CP}. Apply the lower-bound assertion of \cite[Lemma~6.1]{CP} with $h(g)=\sqrt g\log g$. The scale satisfies that as $g\to \infty$,
\[
 \frac{h(g+1)}{h(g)}\longrightarrow1 \quad \textit{$\mathrm{and}$} \quad 
 \frac{\sqrt{\log g}}{h(g)}\longrightarrow0.
\]
Thus the bounded genus shifts and the pinching error in that
lemma do not change the normalized lower limit. For every constant
$c$ below the closed surface lower limit, the cited result gives
$$\liminf\limits_{g\to\infty}\frac{D_{g,n}}{h(g)}\geq c.$$
Letting $c$ increase to that lower limit proves the claim.
\end{proof}

Now we are ready to prove Theorem~\ref{thm:main}.
\begin{proof}[Proof of Theorem~\ref{thm:main}]
By Lemma~\ref{lem:puncture-transfer}, it suffices to treat closed surfaces.
Choose $X\in\M_g^{\geq1}$, which is
nonempty for large $g$ by Proposition~\ref{prop:thickvolume}.
For $R>D_g$, the ball $B_{\M_g}(X,R)$ contains $\M_g^{\geq1}$.
Theorem~\ref{thm:ball-intro} and Proposition~\ref{prop:thickvolume}
give
\[
 2g\log g-Cg
 \leq\log\Vol(\M_g^{\geq1})
 \leq Cg+\sqrt{6A_1g}\,R.
\]
Absorbing the linear error terms gives
\[
 R\geq\frac{2g\log g-Cg}{\sqrt{6A_1g}}.
\]
Let $R\downarrow D_g$ and use $A_1=823\,580$ from
\eqref{eq:numerical-integrated}. We obtain
\[
 D_g\geq\frac{2\sqrt g\log g}{\sqrt{6\cdot823\,580}}-C\sqrt g.
\]
Since $823\,580<4\cdot454^2$, this yields
 \[
 \liminf_{g\to\infty}\frac{\diam(\M_g)}{\sqrt g\log g}
 \geq\frac1{454\sqrt6}
\]
as desired.
\end{proof}


\end{document}